\documentclass[11pt]{article}
\usepackage{amsmath,amssymb,eucal}
\usepackage[dvips]{graphicx}

\newtheorem{thm}{\bf Theorem}[section]

\newtheorem{ex}[thm]{\bf Example}

\numberwithin{equation}{section}

\begin{document}

\title{{\bf 
{\large Seifert surgery on knots via Reidemeister torsion and Casson-Walker-Lescop invariant}}}
\author{
{\small Dedicated to Professor Yukio Matsumoto for his 70th birthday}
\medskip\\
{\normalsize Teruhisa Kadokami, Noriko Maruyama and Tsuyoshi Sakai}}
\date{{\normalsize March 23, 2015}}
\footnotetext[0]{%
}
\maketitle

\begin{abstract}{
For a knot $K$ with $\Delta_K(t)\doteq t^2-3t+1$ in a homology $3$-sphere,
let $M$ be the result of $2/q$-surgery on $K$.
We show that appropriate assumptions on the Reidemeister torsion and
the Casson-Walker-Lescop invariant of the universal abelian covering of $M$ imply $q=\pm 1$, if 
$M$ is a Seifert fibered space.

}\end{abstract}

\section{Introduction}\label{sec:intro}
Dehn surgeries on knots or links have been studied from various points of view
(e.g.\ \cite{Ber, BL, BW, CGLS, Kd1, Kd2, Kd3, KMS, Ma1, Ma2,OS1, OS2, Th, Tr1, Tr2,Wan}).
The first author \cite{Kd1} introduced an idea for applying the Reidemeister torsion to Dehn surgery,
and showed the following:
\begin{thm}\label{th:Seifert}
{\rm (\cite[Theorem 1.4]{Kd2})}
Let $K$ be a knot in a homology $3$-sphere $\Sigma$ such that
the Alexander polynomial of $K$ is $t^2-3t+1$.
The only surgeries on $K$ that may produce a Seifert fibered space with base $S^2$
and with $H_1\ne \{0\}, \mathbb{Z}$ have coefficients $2/q$ and $3/q$,
and produce Seifert fibered space with three singular fibers.
Moreover
(1) if the coefficient is $2/q$, then the set of multiplicities is
$\{2\alpha, 2\beta, 5\}$ where $\gcd(\alpha, \beta)=1$, and
(2) if the coefficient is $3/q$, then the set of multiplicities is
$\{3\alpha, 3\beta, 4\}$ where $\gcd(\alpha, \beta)=1$.
\end{thm}

In this paper, based on Theorem 1.1, we discuss the $2/q$ - Seifert surgery by applying   
the Reidemeister torsion and the Casson-Walker-Lescop invariant in combination simultaneously, and give a sufficient condition to determine the integrality of $2/q$ (Theorem 2.1). The condition is the one suggested by computations for the figure eight knot (Example 2.2).

This paper is actually a continuation of \cite{Kd2}, 
so we follow mainly the notations of \cite{Kd2} and review necessary minimum ones:

\medskip

\noindent
(1) Let $\Sigma$ be a homology $3$-sphere, 
and let $K$ be a knot in $\Sigma$.
Then $\Delta_K(t)$ denotes the Alexander polynomial of $K$, and
$\Sigma(K; p/r)$ denotes the result of $p/r$-surgery on $K$.

\bigskip

\noindent
(2) Let $\zeta_d$ be a primitive $d$-th root of unity.
For an element $\alpha$ of $\mathbb{Q}(\zeta_d)$,
$N_d(\alpha)$ denotes the norm of $\alpha$ associated to
the algebraic extension $\mathbb{Q}(\zeta_d)$ over $\mathbb{Q}$.
Let $f(t)$ be a Laurent polynomial over $\mathbb{Z}$.
We define $|f(t)|_d$ by
$$|f(t)|_d=|N_d(f(\zeta_d))|
=\left| \prod_{i\in (\mathbb{Z}/d\mathbb{Z})^{\times}}
f(\zeta_d^i)\right|.$$
Let $X$ be a homology lens space with $H_1(X)\cong \mathbb{Z}/p\mathbb{Z}$.
Let $d$ be a divisor of $p$.
We define $|X|_d$ by
$$|X|_d=|\Delta_K(t)|_d,$$
where $K$ is a knot in a homology $3$-sphere $\Sigma$ such that
$X=\Sigma(K; p/r)$.
Then $|X|_d$ is a topological invariant of $X$ (Refer to \cite{Kd3} for details).

\bigskip

\noindent
(3) Let $X$ be a closed oriented $3$-manifold.
Then $\lambda(X)$ denotes the Lescop invariant of $X$ (\cite{Le}). Note that $\lambda\left(S^3\right)=0$.

\section{Result}\label{sec:result}

Let $K$ be a knot in a homology $3$-sphere $\Sigma$.
Let $M$ be the result of $2/q$-surgery on $K$: $M=\Sigma(K; 2/q)$.
Let $\pi : X\to M$ be the universal abelian covering of $M$
(i.e.\ the covering associated to $\mathrm{Ker}(\pi_1(M)\to H_1(M))$).
Since $H_1(M)\cong \mathbb{Z}/2\mathbb{Z}$, $\pi$ is the $2$-fold unbranched covering.

We then define $\lambda_q(K)$ by the following formula:
$$\lambda_q(K):=\lambda(X).$$
It is obvious that $\lambda_q(K)$ is a knot invariant of $K$. 
We also define $|K|_{(q,d)}$ by the following formula, if $|X|_d$ is defined:
$$|K|_{(q,d)}:=|X|_d.$$
It is also obvious that $|K|_{(q,d)}$ is a knot invariant of $K$.

We then have the following.

\begin{thm}\label{th:main}
Let $K$ be a knot in a homology 3-sphere $\Sigma$. We assume the following.
 
\medskip
\noindent
{\bf (2.1)} $\lambda(\Sigma)=0$,

\medskip
\noindent
{\bf (2.2)} $\Delta_K(t)\doteq t^2-3t+1$,

\medskip
\noindent
{\bf (2.3)} $|q|\not =1$,
 
\medskip
\noindent
{\bf (2.4)} $\sqrt{|K|_{(q,5)}}\geq 4\{\lambda_q(K)\}^2-1$.
\medskip

Then $M=\Sigma(K; 2/q)$ is not a Seifert fibered space.
\end{thm}

The assumption {\bf (2.2)} implies $H_1(X)\cong \mathbb{Z}/5\mathbb{Z}$ as shown in \S 3, 
hence $|X|_5$ is defined. The assumption {\bf (2.4)} means $\sqrt{|X|_5}\geq 4\{\lambda(X)\}^2-1$. As noticed in the introduction, this inequality is
suggested by computations in the following example, and so is the assumption {\bf (2.1)}.

\begin{ex}\label{ex:fig8}
{\rm Let $K$ be the figure eight knot. Then $\lambda_q(K)=-q$ and $|K|_{(q,5)}=(5q^2-1)^2$. Hence {\bf (2.4)} holds for every $q$.}
\end{ex}

\section{Proof of Theorem \ref{th:main}}\label{sec:proof}

Let $\Sigma_2$ be the double branched covering space of $\Sigma$ branched along $K$,
and $\overline{K}$ the lifted knot of $K$ in $\Sigma_2$.
Since $|\Delta_K(-1)|=5$, we have $H_1(\Sigma_2)\cong \mathbb{Z}/5\mathbb{Z}$.
Since $\overline{K}$ is null-homologous in $\Sigma_2$,
and $X$ is the result of $1/q$-surgery on $\overline{K}$,
we have $H_1(X)\cong \mathbb{Z}/5\mathbb{Z}$.

\medskip

We suppose that $M$ is a Seifert fibered space.
According to Theorem \ref{th:Seifert}, we may assume that
$M$ has a framed link presentation as in Figure 1, 
where $1\le \alpha<\beta$ and $\gcd(\alpha, \beta)=1$.

\begin{figure}[htbp]
\begin{center}
\includegraphics[scale=0.6]{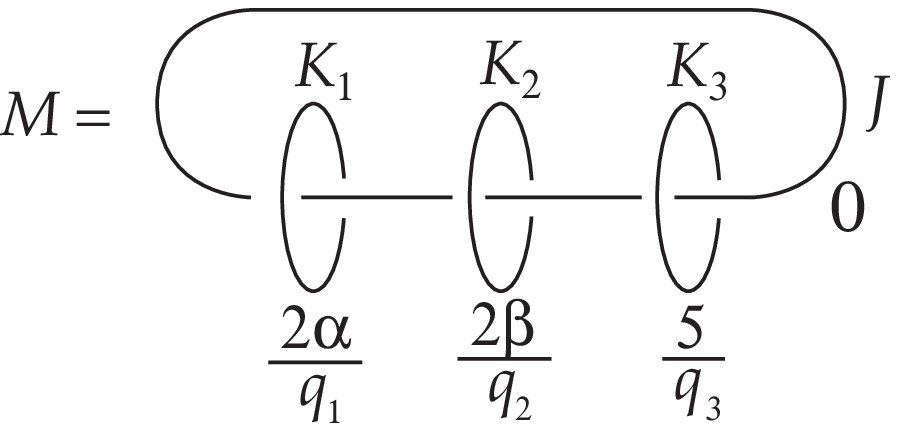}
\label{M1}
\caption{A framed link presentation of $M=\Sigma (K; 2/q)$}
\end{center}
\end{figure}

\medskip
\noindent
In fact, since $H_1(M)\cong \mathbb{Z}/2\mathbb{Z}$, 
the base surface of $M$ has genus $0$ and is $S^2$
(i.e.\ the projective plane $\mathbb{P}^2$ is ruled out because
the order of $H_1(M)$ is not divisible by $4$).
By {\bf (2.2)}, $\Delta_K(t)\doteq t^2-3t+1$.
Hence by Theorem \ref{th:Seifert},
we have the presentation in Figure 1, 
where $1 \leq \alpha \leq \beta$ and $\gcd(\alpha, \beta)=1$. Then the exceptional case that $\alpha=\beta=1$ is removed as follows:
The order of $H_1(M)$ is equal to
$$\left|10q_1+10q_2+4q_3\right|$$
if $\alpha=\beta=1$. Since $H_1(X)\cong \mathbb{Z}/2\mathbb{Z}$, we have
$$5q_1+5q_2+2q_3=\pm 1$$
Since $q_1$ and $q_2$ are odd, this is impossible. Thus we have the desired presentation.  

\medskip

Then, on the universal abelian covering $X$ of $M$, we see

\begin{center}
$(\ast)$ : $X$ has a framed link presentation as in Figure 2.
\end{center}

\begin{figure}[htbp]
\begin{center}
\includegraphics[scale=0.6]{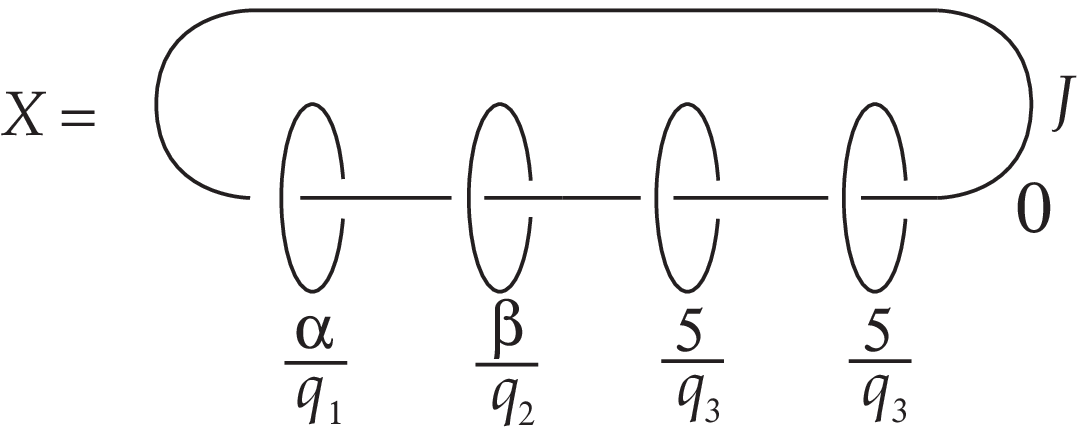}
\label{X}
\caption{A framed link presentation of $X$}
\end{center}
\end{figure}

Assuming $(\ast)$, we proceed with the proof (We will give a proof of $(\ast)$ in the appendix, 
which is essentially owing to \cite{Se}).

\medskip

By $(\ast)$ and \cite[Theorem 1.2 (3)]{Kd2},
we have $|X|_5=(\alpha \beta)^4$.
Hence by {\bf (2.4)}, we have
\begin{equation}\label{eq:Pell}
(\alpha \beta)^2\geq 4\{\lambda(X)\}^2-1.
\end{equation}
Since $4\{\lambda(X)\}^2\leq (\alpha \beta)^2+1<(\alpha \beta+1)^2$, we have
\begin{equation}\label{eq:ineq1}
{\displaystyle
|\lambda(X)|<\frac{\alpha \beta+1}{2}.}
\end{equation}

We now consider $e$ defined as follows:
$$e:=\frac{q_1}{\alpha}+\frac{q_2}{\beta}+\frac{q_3}{5}+\frac{q_3}{5}.$$

According to the sign of $e$, we treat two cases separetely: We first consider the case $e>0$.
Then the order of $H_1(X)$ is $25\alpha \beta e$
by $(\ast)$ and \cite{Or, Sv}.
Since $H_1(X)\cong \mathbb{Z}/5\mathbb{Z}$, $25\alpha \beta e=5$, and
$e=1/(5\alpha \beta)$.
Hence by $(\ast)$ and \cite[Proposition 6.1.1]{Le}, we have
\begin{equation}\label{eq:neweq}
{\displaystyle
\lambda(X)=(-2)\alpha \beta
+\frac{25\beta}{24\alpha}
+\frac{25\alpha}{24\beta}
+\frac{1}{24\alpha \beta}-\frac 58-\frac 52 S
}
\end{equation}
where $S=s(q_1, \alpha)+s(q_2, \beta)+2s(q_3, 5)$
and $s(\cdot, \cdot)$ denotes the Dedekind sum ([RG]).

By (\ref{eq:ineq1}), we have
$$-\frac{\alpha \beta+1}{2}<\lambda(X).$$
Hence 
$$-\frac{\alpha \beta+1}{2}<
(-2)\alpha \beta+\frac{25\beta}{24\alpha}+\frac{25\alpha}{24\beta}
+\frac{1}{24\alpha \beta}-\frac 58+\frac 52 |S|.$$
Consequently
we have
\begin{equation}\label{eq:ineq2}
{\displaystyle
\frac 32\alpha \beta
<-\frac 18+\frac{25}{24\alpha}\beta+\frac{25}{24}\left( \frac{\alpha}{\beta}\right)
+\frac{1}{24\alpha \beta}+\frac 52 |S|.
}
\end{equation}

We show that $\alpha \ge 2$ implies a contradiction: 
Suppose that $\alpha \ge 2$.
Since $\alpha<\beta$, we have $\beta \ge 3$ and $\alpha/\beta<1$.
Hence 
$$3\beta<-\frac 18+\frac{25}{24\cdot 2}\beta+\frac{25}{24}+\frac{1}{24\cdot 2\cdot 3}
+\frac 52 |S|.$$
In general, $|s(q, p)|\le p/12$ holds.
In fact, by \cite{BL} and \cite{RG},
$$|s(q, p)|\le s(1, p)=\frac{(p-1)(p-2)}{12p}\le \frac{p^2}{12p}=\frac{p}{12}.$$
Since $|s(q_1,\alpha)| \leq \frac \alpha {12} < \frac \beta {12}$, $|s(q_2,\beta)| \leq \frac \beta {12}$, 
and $|s(q_3,5)| \leq \frac 15$, we have  
$$|S|\le |s(q_1, \alpha)|+|s(q_2, \beta)|+2|s(q_3, 5)|
\le \frac{\beta}{6}+\frac{2}{5}.$$
Hence
\begin{eqnarray*}
3\beta & < & -\frac 18+\frac{25}{48}\beta+\frac{25}{24}+\frac{1}{144}
+\frac 52\left( \frac{\beta}{6}+\frac 25\right) \bigskip\\
& = & \left( \frac 78+\frac{1}{144}\right)+\frac{25}{24}
+\frac{45}{48}\beta \bigskip\\
& < & 1+1+2+\beta
\end{eqnarray*}
implies $\beta<2$.
This contradicts $2\le \alpha<\beta$.

\medskip

We next show that $\alpha=1$ implies $\beta<6$: 
Suppose that $\alpha=1$, then $\beta \ge 2$. Substituing $\alpha=1$ in 
(\ref{eq:ineq2}), we have
$$\frac 32\beta<-\frac 18+\frac{25}{24}\beta+\frac{25}{24\beta}
+\frac{1}{24\beta}+\frac 52 |S|$$
where $S=s(q_2, \beta)+2s(q_3, 5)$
(since $s(q_1, 1)=0$).
By using $\beta \ge 2$,
$$\frac 32\beta<-\frac 18+\frac{25}{24}\beta+\frac{25}{24\cdot 2}
+\frac{1}{24\cdot 2}+\frac 52 |S|.$$
By applying $|s(q, p)|\le \frac p{12}$ and $|s(q_3,5)| \leq \frac 15$,
$$|S|\le |s(q_2, \beta)|+2|s(q_3, 5)|\le \frac{\beta}{12}+\frac 25,$$
and hence
\begin{eqnarray*}
\frac 32 \beta & < & -\frac 18+\frac{25}{24}\beta+\frac{26}{48}
+\frac 52 \left( \frac{\beta}{12}+\frac 25\right)
\bigskip\\
& = & \frac{17}{12}+\frac 54\beta.
\end{eqnarray*}
Thus we have $\beta<6$.

Since $\alpha=1$, $e=\frac 1{5\beta}$. Hence
$$\frac{q_1}{1}+\frac{q_2}{\beta}+\frac{q_3}{5}+\frac{q_3}{5}=\frac 1{5\beta},$$
and hence we have the following equation
\begin{equation}\label{eq:neweq2}
(5\beta)q_1+5q_2+(2\beta)q_3=1.
\end{equation}
Since $q_1$ and $q_2$ are odd (see Figure 1), 
$\beta$ must be even by (\ref{eq:neweq2}). 
Hence $\beta=2$ or $4$.

Suppose that $\beta =2$, then by (\ref{eq:neweq2}),
$$10q_1+5q_2+4q_3=1.$$
Hence $q_2\equiv -1 \;(\mathrm{mod}\ \! 4)$ and $q_3\equiv -1 \;(\mathrm{mod}\ \! 5)$. 
By using [Le, Proposition 6.1.1] again,
we have $\lambda(M)=-1$.
On the other hand, $\lambda(M)=-q$ by the assumptions {\bf (2.1)} and {\bf (2.2)}. 
Hence $q=1$. This contradicts {\bf (2.3)}.
 
Suppose that $\beta =4$, then by (\ref{eq:neweq2}),
$$20q_1+5q_2+8q_3=1.$$
Hence $q_2\equiv 1 \;(\mathrm{mod}\ \! 4)$ and $q_3\equiv 2 \;(\mathrm{mod}\ \! 5)$. 
Hence by (\ref{eq:neweq}),
we have $\lambda(X)=-\frac 92$. This contradicts (\ref{eq:Pell}), 
and ends the proof in the case $e>0$.

\medskip

We finally consider the case $e<0$.
Then ${\displaystyle e=-\frac{1}{5\alpha \beta}}$.
By $(\ast)$ and \cite[Proposition 6.1.1]{Le}, we have
$$\lambda(X)=-\left\{(-2)\alpha \beta
+\frac{25\beta}{24\alpha}+\frac{25\alpha}{24\beta}
+\frac{1}{24\alpha \beta}-\frac 58+\frac 52 S
\right\}.$$
Remaining part of the proof is similar to that in the case $e>0$.

\medskip

This completes the proof of Theorem 2.1.\hfill
$\Box$

\section{A formula for the Casson-Walker-Lescop invariant
of the result of surgery on $2$-bridge link}\label{sec:CWL}

\begin{figure}[htbp]
\begin{center}
\includegraphics[scale=0.5]{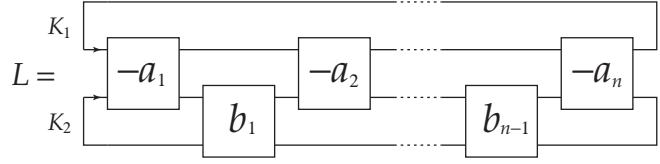}
\label{twobrl}
\caption{$2$-bridge link $L=K_1\cup K_2=D(a_1, b_1, a_2, b_2, \ldots, a_{n-1}, b_{n-1}, a_n)$}
\end{center}
\end{figure}
Let $L=K_1\cup K_2=D(a_1, b_1, a_2, b_2, \ldots, a_{n-1}, b_{n-1}, a_n)$ be 
the oriented $2$-bridge link as in Figure 3, where $(a_1, b_1, \ldots, a_n)$
is a sequence of $(2n-1)$ integers, and a frame box labeled an integer $c$
denotes a $|c|$ positive (resp. negative) full twists
of two horizontal strands if $c>0$ (resp. $c<0$).

\medskip

Let $X$ be the result of surgery on $L$ with coefficients $p_1/q_1$ and $p_2/q_2$
being $q_1$ and $q_2$ positive:
$$X=S^3(L; p_1/q_1, p_2/q_2),\quad
q_1>0,\ q_2>0.$$
Let $E$ denote the linking matrix associated to the framed link presentation of $X$ defined by
$$E=\left( 
\begin{array}{cc}
p_1/q_1 & \ell \\
\ell & p_2/q_2
\end{array}\right)$$
where $\ell$ is the linking number of $L$.
Note that $\ell=-\sum_{i=1}^na_i$.
Let $\mathrm{tr}(E)$, $\sigma(E)$ and $b_-(E)$ denote
the trace, the signature and the number of negative eigenvalues of $E$, respectively.
Then the following formula holds:

\begin{thm}\label{th:CW}
{\rm (\cite[Lemma 2.1 and Proposition 3.1]{Ma2})}
$$\lambda(X)
=(-1)^{b_-(E)}q_1q_2
\left(
\frac{p_2}{q_2}[K_1]+\frac{p_1}{q_1}[K_2]+[L]
\right)
+|p|\left(
\frac{\sigma(E)}{8}+\frac{s(p_1, q_1)}{2}+\frac{s(p_2, q_2)}{2}
\right)$$
where 
$[L]=\sum_{k=1}^{n-1}b_k(a_1+\cdots+a_k)(a_{k+1}+\cdots+a_n)
-\frac{\ell(\ell^2-1)}{12}+\frac{\ell^2}{12}\mathrm{tr}(E)$, 
$[K_i]=-\frac{p_i^2+q_i^2+1}{24q_i^2}$ for $i=1, 2$, and
$|p|=q_1q_2|\det(E)|$.
\end{thm}
We apply this theorem in the next section.

\section{Proof of Example \ref{ex:fig8}}\label{sec:example}
It is well-known that the figure eight knot $4_1$ satisfies {\bf (2.2)} of Theorem 2.1.

\medskip

We now consider the oriented link $L=K_1\cup K_2$ shown in Figure 4.

\medskip

\begin{figure}[htbp]
\begin{center}
\includegraphics[scale=0.5]{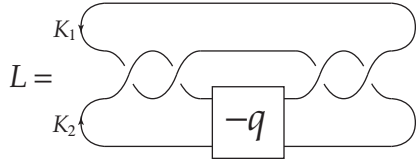}
\label{L}
\caption{$L=K_1\cup K_2=D(1, -q, 1)$}
\end{center}
\end{figure}

\noindent 
Then $L=D(1, -q, 1)$ in the notation of Section \ref{sec:CWL}.

\medskip

As is well-known \cite{Ro}, $M=S^3(4_1; 2/q)$ is presented as shown in Figure 5,
\begin{figure}[htbp]
\begin{center}
\includegraphics[scale=0.45]{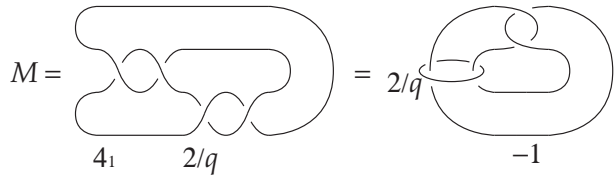}
\label{M2}
\caption{$M=S^3(4_1; 2/q)$}
\end{center}
\end{figure}
and hence its two fold unbranched covering $X$ is presented as $S^3(L; -3, -3)$  shown in Figure 6.
\begin{figure}[htbp]
\begin{center}
\includegraphics[scale=0.45]{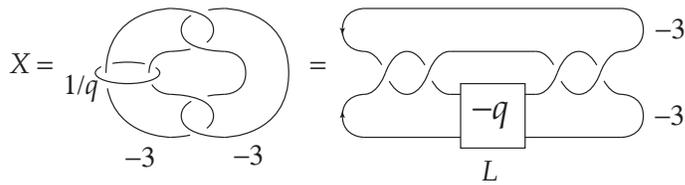}
\label{X}
\caption{$X=S^3(L; -3, -3)$}
\end{center}
\end{figure}

\medskip
We compute $\lambda(X)$: Let $E$ be the linking matrix associated to the above framed link presentation of $X$:
$E=\left(
\begin{array}{cc}
-3 & -2\\
-2 & -3
\end{array}
\right)$.
Since $\mathrm{tr}(E)=-6$, $b_-(E)=2$, $\sigma(E)=-2$
, $[L]=-q-\frac 32$ and $[K_i]=-\frac{11}{24}$ $(i=1,2)$, 
we have $\lambda(X)=-q$ by Theorem \ref{th:CW}.

\medskip

We finally compute $|X|_5$: 
According to  \cite[Theorem 1]{Kn}, 
the $2$-variable Alexander polynomial $\Delta_L(t_1, t_2)$ of $L$ is given by
\begin{equation}\label{eq:Alx}
\Delta_L(t_1, t_2)\doteq -q(t_1-1)(t_2-1)+t_1t_2+1.
\end{equation}

\noindent 
Further, let $T_i$ $(i=1, 2)$ be the representing element of a meridian of $K_i$ in $H_1(X)$.
Then by the linking matrix $E$, we see
$$H_1(X)=\langle T_1, T_2\ |\ T_1^{-3}T_2^{-2}=T_1^{-2}T_2^{-3}=1\rangle
=\langle T\ |\ T^5=1\rangle \cong \mathbb{Z}/5\mathbb{Z},$$
where we regard $T=T_1=T_2$.

Let $\zeta=\zeta_5$ be a primitive $5$-th root of unity, and
$\varphi : \mathbb{Z}[H_1(X)]\to \mathbb{Q}(\zeta)$
a ring homomorphism defined by $\varphi(T_1)=\varphi(T_2)=\zeta$.
Then the Reidemeister torsion of $X$ associated to $\varphi$, denoted by $\tau^{\varphi}(X)$, 
is defined (cf.\ \cite{Tr1, Tr2}), and we have
\begin{equation}\label{eq:tor1}
\tau^{\varphi}(X)\doteq
\{(1-q)(\zeta-1)^2+2\zeta\}(\zeta-1)^{-2}
\end{equation}
by (\ref{eq:Alx}) and \cite[Lemma 2.5 (1)]{Kd2}.

In addition, 
suppose that $X=\Sigma'(K'; 5/q')$ for a knot $K'$ in a homology $3$-sphere $\Sigma'$.
Then we have
\begin{equation}\label{eq:tor2}
\tau^{\varphi}(X)\doteq
\Delta_{K'}(\zeta')(\zeta'-1)^{-1}(\zeta'^{{\bar q}'}-1)^{-1}
\end{equation}
by \cite[Lemma 2.6]{Kd2},
where $\zeta'$ is a primitive $5$-th root of unity, and $q'{\bar q}'\equiv 1\ (\mathrm{mod}\ \! p)$.

By comparing (\ref{eq:tor1}) and (\ref{eq:tor2}), we have
$$\begin{array}{ll}
|X|_5=|\Delta_{K'}(t)|_5
& {\displaystyle
=\left|\prod_{i\in (\mathbb{Z}/5\mathbb{Z})^{\times}}
\{(1-q)(\zeta^i-1)^2+2\zeta^i\}\right|}
\bigskip\\
& =\{(1-q)(\zeta+\zeta^{-1})+2q\}^2\{(1-q)(\zeta^2+\zeta^{-2})+2q\}^2
\bigskip\\
& =(5q^2-1)^2.
\end{array}$$
Since $\lambda(X)=-q$, we have
\begin{eqnarray}
\sqrt{|X|_5}&=&5\{\lambda(X)\}^2-1 \nonumber \\
&\geq& 4\{\lambda(X)\}^2-1. \nonumber
\end{eqnarray}
Therefore the figure eight knot satisfies {\bf (2.4)}.\hfill{$\Box$}

\section{Concluding Remarks}\label{sec:remark}

\noindent
(1) W.~Thurston \cite{Th} determined exceptional surgeries of the figure eight knot
in terms of hyperbolic geometry, and
M.~Brittenham and Y.~Wu \cite{BW} determined exceptional surgeries 
of $2$-bridge knots by using lamination theory.

\bigskip

\noindent
(2) As in Theorem \ref{th:Seifert}, 
the abelian Reidemeister torsion of $M$ dominates the numerator $p$ of
Seifert surgery coefficient $p/q$ on a knot.
On the other hand, as in Theorem \ref{th:main} and Example \ref{ex:fig8},
the meta-abelian Reidemeister torsion $\tau^{\varphi}(X)$ of $M$
and the Casson-Walker-Lescop invariant $\lambda(X)$ in combination dominate the denominator $q$ of $p/q$.

\bigskip

\noindent
(3) M.~Marcolli and B.~Wang \cite{MW}, and L.~Nicolaescu \cite{Ni} showed
that the Seiberg-Witten invariant for a rational homology $3$-sphere
is decomposed into the Reidemeister-Turaev torsion part
and the Casson-Walker invariant part.

\medskip

Hence it would be worth asking how one can prove directly 
Theorem \ref{th:main} or the like
by applying the Seiberg-Witten invariant.

\bigskip

\noindent
{\bf Acknowledgement}\ 
The authors would like to thank to 
Kazuhiro Ichihara, Kimihiko Motegi, Makoto Sakuma and
Yasuyoshi Tsutsumi for giving them useful comments.

\newpage

\noindent
{\bf \Large Appendix: Proof of $(\ast)$ in Section 3}\label{sec:appendix}

\bigskip

\begin{figure}[htbp]
\begin{center}
\includegraphics[scale=0.5]{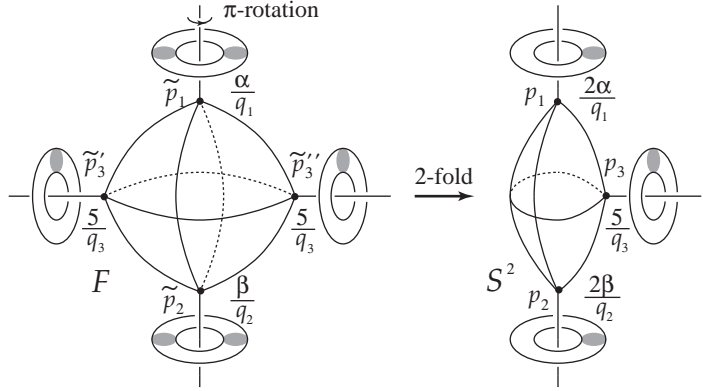}
\label{branch}
\caption{Illustration for proof}
\end{center}
\end{figure}

Figure 7 illustrates the proof given below:
Let $N$ be the result of $0$-surgery along $J$, and
we consider $K_1\cup K_2\cup K_3$ as a link in $N$ (Figure 8 (a)).
\begin{figure}[htbp]
\begin{center}
\includegraphics[scale=0.6]{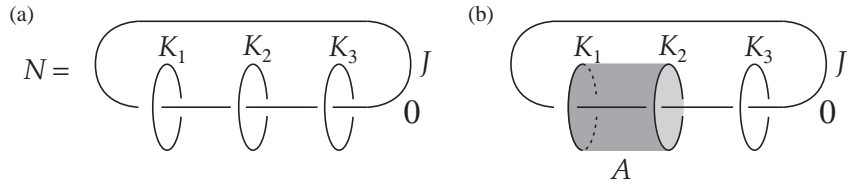}
\label{N}
\caption{$K_1\cup K_2\cup K_3$ in $N=S^2\times S^1$, and
an annulus $A$ bounded by $K_1\cup K_2$}
\end{center}
\end{figure}
Then $N=S^2\times S^1$, and we may assume that
$K_1$, $K_2$ and $K_3$ are regular fibers of $N$; 
in other words we may assume as follows:
$$K_1=\{p_1\}\times S^1,\ K_2=\{p_2\}\times S^1,\ K_3=\{p_3\}\times S^1,
\quad p_1, p_2, p_3\in S^2.$$
Note that we can choose a regular fiber of $N$ as a preferred longitude of $K_i$ in Figure 1.
In the following, we always choose a regular fiber as a preferred longitude of $K_i$.

\medskip

Let $A$ be an annulus as illustrateted in Figure 8 (b),
which is a Seifert surface for $K_1\cup K_2$ in $N$.
Let $Y$ be the $2$-fold branched covering of $N$ branched along $K_1\cup K_2$
constructed by cut-open and copy-paste along $A$.
Then $Y$ is a Seifert fibered space.
Let $F$ be the base surface of $Y$ .
Then $F$ is the $2$-fold branched covering of $S^2$ with $\{p_1, p_2\}$
as the branch set.
Hence $F=S^2$, and $Y=F\times S^1=S^2\times S^1$. 

\medskip

Let $\widetilde{p}_1$, $\widetilde{p}_2$, $\{\widetilde{p}'_3, \widetilde{p}''_3\}$
be the inverse images of $p_1$, $p_2$, $p_3$ with respect to
$F\to S^2$ respectively, and set as follows:
$$\widetilde{K}_1=\{\widetilde{p}_1\}\times S^1,\ 
\widetilde{K}_2=\{\widetilde{p}_2\}\times S^1,\ 
\widetilde{K}'_3=\{\widetilde{p}'_3\}\times S^1,\ 
\widetilde{K}''_3=\{\widetilde{p}''_3\}\times S^1.$$
Then $\widetilde{K}_1$, $\widetilde{K}_2$, $\widetilde{K}'_3\cup \widetilde{K}''_3$
are the inverse images of $K_1$, $K_2$, $K_3$ with respect to
$Y\to N$ respectively.
Note that 
$$\widetilde{K}_1\cup \widetilde{K}_2\cup \widetilde{K}'_3\cup \widetilde{K}''_3
\subset Y=S^2\times S^1$$
is viewed as in Figure 9.

\begin{figure}[htbp]
\begin{center}
\includegraphics[scale=0.6]{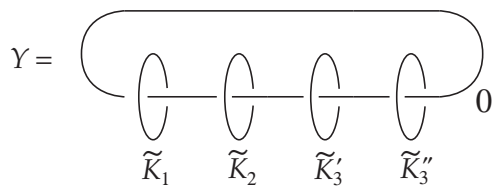}
\label{Y}
\caption{$2$-fold branched covering $Y\to N$ along $\partial A$, and
the lifted link of $K_1\cup K_2\cup K_3$}
\end{center}
\end{figure}

By Figure 1 and Figure 8 (a), 
$M$ is obtained from $N$ by surgery on $K_1\cup K_2\cup K_3$
with coefficients $2\alpha/q_1$, $2\beta/q_2$, $5/q_3$.
We have assumed that the preferred longitude of $K_i$ in Figure 1
is a regular fiber of $N$, and hence its lift to $Y$ is also a regular fiber
(since each lift of a regular fiber is a regular fiber).
Hence $X$ is obtained from $Y$ by surgery on 
$\widetilde{K}_1\cup \widetilde{K}_2\cup \widetilde{K}'_3\cup \widetilde{K}''_3$
with coefficients $\alpha/q_1$, $\beta/q_2$, $5/q_3$, $5/q_3$.
This means $(\ast)$.\hfill{$\Box$}

{\small
\par
Teruhisa Kadokami\par
Department of Mathematics, East China Normal University,\par
Dongchuan-lu 500, Shanghai, 200241, China\par
{\tt mshj@math.ecnu.edu.cn, kadokami2007@yahoo.co.jp}\par

\medskip

Noriko Maruyama\par
Musashino Art University,\par 
Ogawa 1-736, Kodaira, Tokyo 187-8505, Japan \par 
{\tt maruyama@musabi.ac.jp} \par

\medskip

Tsuyoshi Sakai\par
Department of Mathematics, Nihon University,\par
3-25-40, Sakurajosui, Setagaya-ku, Tokyo 156-8550, Japan \par
}


\bigskip

{\footnotesize
 \begin{thebibliography}{999}
%
\newcommand{\bysame}{%
       \leavevmode\hbox to 3em{\hrulefill}\,}
%
\bibitem[Ber]{Ber}J.~Berge,
{\rm Some knots with surgeries yielding lens spaces},
(Unpublished manuscript, 1990).
%
\bibitem[BL]{BL}S.~Boyer and D.~Lines,
{\it Surgery formulae for Casson's invariant and extensions to homology lens spaces},
J. Reine Angew. Math., {\bf 45}\ (1990), 181--220.
%
\bibitem[BW]{BW}M.~Brittenham and Y.~Wu, 
{\it The classification of exceptional Dehn surgeries on $2$-bridge knots},
Comm. Anal. Geom. {\bf 9}\ (2001), 97--113.
%
\bibitem[CGLS]{CGLS}M.~Culler, M.~Gordon, J.~Luecke and P.~Shalen, 
{\it Dehn surgery on knots},
Ann. of Math., {\bf 125} (1987), 237--300.
%
\bibitem[Kd1]{Kd1}T.~Kadokami,
{\it Reidemeister torsion and lens surgeries on knots in homology $3$-spheres I},
Osaka J. Math., {\bf 43}, no.4\ (2006), 823--837.
%
\bibitem[Kd2]{Kd2}T.~Kadokami,
{\it Reidemeister torsion of Seifert fibered homology lens spaces and Dehn surgery},
Algebr. Geom. Topol., {\bf 7}\ (2007), 1509--1529.
%
\bibitem[Kd3]{Kd3}T.~Kadokami,
{\it Reidemeister torsion and lens surgeries on knots in homology $3$-spheres II},
Top. Appl., {\bf 155}, no.15\ (2008), 1699--1707.
%
\bibitem[KMS]{KMS}T.~Kadokami, N.Maruyama and M.~Shimozawa, {\it Lens surgeries along the $n$-twisted Whitehead link}, Kyunpook Math. J., {\bf 52} (2012), 245--264.
%
\bibitem[Kn]{Kn}T.~Kanenobu, 
{\it Alexander polynomials of two-bridge links},
J. Austral. Math. Soc. (Ser. A), {\bf 36}\ (1984), 59--68.
%
\bibitem[Le]{Le}C.~Lescop, 
{\it Global surgery formula for the Casson-Walker invariant},
Ann. of Math. Studies, Princeton Univ. Press., {\bf 140} (1996).
%
\bibitem[MW]{MW}M.~Marcolli and B.~Wang, 
{\it Seiberg-Witten and Casson-Walker invariants for rational homology $3$-spheres},
Geom. Dedicata, {\bf 91} (2002), 45--58.
%
\bibitem[Ma1]{Ma1}N.~Maruyama, 
{\it On Dehn surgery along a certain family of knots},
Jour. of Tsuda College, {\bf 19}\ (1987),  261--280.
%
\bibitem[Ma2]{Ma2}N.~Maruyama, 
{\it The CWL invariant and surgeries along $2$-component links},
Jour. of Musashino Art University, {\bf 42}\ (2011),  39--51.
%
\bibitem[Ni]{Ni}L.~Nicolaescu, 
{\it Seiberg-Witten invariants of rational homology $3$-spheres},
Commun. Contemp. Math. {\bf 6}, No. 6\ (2004), 833--866.
%
\bibitem[Or]{Or}P.~Orlik, 
{\it Seifert manifolds},
Lecture Notes in Math. {\bf 291} (1972), Springer-Verlag.
%
\bibitem[OS1]{OS1}P.~Ozsv\'ath and Z.~Szab\'o, 
{\it On Heegaard Floer homology and Seifert fibered surgeries},
Geom. Top., {\bf 7} (2004), 181--203.

\bibitem[OS2]{OS2}P.~Ozsv\'ath and Z.~Szab\'o, 
{\it On knot Floer homology and lens space surgeries},
Topology, {\bf 44} (2005), 1281--1300.
%
\bibitem[RG]{RG}H.~Rademacher and E.~Grosswald, 
{\it Dedekind sums},
The Carus Mathematical Monograph, {\bf 16}\ (1972).
%
\bibitem[Ro]{Ro}D.~Rolfsen, 
{\it Knots and Links},
Publish or Perish, Inc. (1976);
AMS Chelsea Publishing (1976).
%
\bibitem[Sv]{Sv}N.~Saveliev, 
{\it Invariants for homology 3-spheres},
Encyclopaedia of Mathematical Sciences, {\bf 140},
Low-Dimensional Topology I, Springer-Verlag Berlin,\ (2002).
%
\bibitem[Se]{Se}H.~Seifert, 
{\it Topologie Dreidimensionaler Gefaserter R\"aume} (German),
Acta. Math., {\bf 60}, no.1\ (1933), 147--238.
%
\bibitem[Th]{Th}W.~P.~Thurston, 
{\it The Geometry and Topology of Three-Manifolds},
Electronic version 1.1 (www.msri.org/publications/books/gt3m/)\ (2002).
%
\bibitem[Tr1]{Tr1}V.~G.~Turaev, 
{\it Reidemeister torsion in knot theory},
Russian Math. Surveys, {\bf 41-1}\ (1986), 119--182.
%
\bibitem[Tr2]{Tr2}V.~G.~Turaev, 
{\it Introduction to Combinatorial Torsions},
Birkh\"auser Verlag\ (2001).
%
\bibitem[Wan]{Wan}J.~Wang,
{\it Cosmetic surgeries on genus one knots},
Algebr. Geom. Topol., {\bf 6}\ (2006), 1491--1517.
%
\end{thebibliography}}
\end{document}